\def\mapright#1{\smash{\mathop{\longrightarrow}\limits^{#1}}}
\def\mapleft#1{\smash{\mathop{\longleftarrow}\limits^{#1}}}

\def\mapup#1{\uparrow\rlap{$\vcenter{\hbox{$\scriptstyle#1$}}$}}
\def\mapupleft#1{\uparrow\llap
                   {$\vcenter{\hbox{$\scriptstyle#1$}}\;\;$}}

\def\simleftarrow{\smash{\mathop{\leftarrow}\limits^{\sim}}}

\def\simmaprightarrow#1{\smash{\mathop{\longrightarrow}\limits^{\sim}_{#1}}}
\def\simmapleftarrow#1{\smash{\mathop{\longleftarrow}\limits^{\sim}_{#1}}}

\def\LKIII{\Lambda_{\text{K3}}}
\def\PE{\mathcal P(T,\Lambda)}
\def\PEj{\mathcal P_{j}(T,\Lambda)}
\def\GPE{\mathcal P^{G\text{-eq}}(T,\Lambda)}
\def\GPEj{\mathcal P^{G\text{-eq}}_j(T,\Lambda)}

\def\GEj{\mathcal E^{G\text{-eq}}_j}

\def\PEKIII{\mathcal P(T,\Lambda_{\text{K3}})}
\def\GPEKIII{\mathcal P^{G\text{-eq}}(T,\Lambda_{\text{K3}})}

\documentclass{amsart}

\usepackage{amscd, amsmath, amsthm, amssymb}

\newtheorem{theorem}{Theorem}[section]
\newtheorem{lemma}[theorem]{Lemma}
\newtheorem{proposition}[theorem]{Proposition}

\newtheorem{corollary}[theorem]{Corollary}

\theoremstyle{definition}     
\newtheorem{definition}[theorem]{Definition}
\newtheorem{definitionThm}[theorem]{Definition-Theorem}

\theoremstyle{remark}

\numberwithin{equation}{section}

\newcommand{\simrightarrow}{\smash{\mathop{\rightarrow}\limits^{\sim}}}

\begin{document}

\title[Fourier-Mukai number of a K3 surface]{
Fourier-Mukai number of a K3 surface}

\author[S. Hosono]{Shinobu Hosono}

\address{
Graduate School of Mathematical Sciences,
University of Tokyo, Komaba Meguro-ku,
Tokyo 153-8914, Japan
}
\email{hosono@ms.u-tokyo.ac.jp}

\author[B.H. Lian]{Bong H. Lian}

\address{
Department of Mathematics, Brandeis University,
Waltham, MA 02154, U.S.A.}
\email{lian@brandeis.edu}

\author[K. Oguiso]{Keiji Oguiso}

\address{
Graduate School of Mathematical Sciences,
University of Tokyo, Komaba Meguro-ku,
Tokyo 153-8914, Japan
}
\email{oguiso@ms.u-tokyo.ac.jp}

\author[S.-T. Yau]{Shing-Tung Yau}

\address{
Department of Mathematics, Harvard University,
Cambridge, MA 02138, U.S.A.}
\email{yau@math.harvard.edu}


\begin{abstract}
We shall give a Counting Formula for the number
of Fourier-Mukai partners of a K3 surface
and consider three applications.
\end{abstract}

\maketitle


\setcounter{section}{-1}
\section{ Introduction}

In his papers [Mu 1,2,3], Mukai discovered the fundamental importance
of manifolds with equivalent bounded derived categories of coherent sheaves.
Such manifolds are now known as Fourier-Mukai (FM) partners, which
have been the focus of attention in several contexts of mathematics.

The aim of this note is to derive an explicit counting formula for
the FM number, i.e. the cardinality
of the set $FM(X)$ of FM partners,
for a projective K3 surface $X$ in terms of
the N\'eron-Severi
lattice $NS(X)$ and the Hodge structure $(T(X), \mathbf C \omega_{X})$
of the transcendental lattice $T(X)$ (Counting Formula 2.3), and to give
then three applications (Corollary 2.7, Theorems 3.3 and 4.1).

The Counting Formula 2.3 is inspired by an earlier work of Bridgeland,
Maciocia [BM]. Our proof here is based on the fundamental
result of Mukai [Mu3] and Orlov [Or]  (see also [BM], [HLOY1] and Theorem 2.2
for summary), and the theory of primitive embeddings of lattices
due to Nikulin [Ni1]. What we actually need here
is a sort of equivariant version of the theory in [Ni1].
Namely, we introduce the notion of $G$-equivalence classes of primitive
embeddings of an even non-degenerate lattice $T$ into an
even unimodular indefinite lattice $\Lambda$, where $G$ is a prescribed
subgroup of the orthogonal group $O(T)$ (Definition 1.1), and
we find a formula for the cardinality of the set $\GPE$ of
$G$-equivalence classes of the primitive
embeddings $\iota : T \rightarrow \Lambda$ (Theorem 1.4).

\vskip0.3cm
We give a precise formula for $\vert \GPE \vert$ in \S 1.
Such a natural and useful formula like Theorem 1.4
is clearly not unexpected.
However, we are unable to find this general formula anywhere
in the literature, and so we give a proof in Appendix A.
We remark that in the special case
$G = \{ \text{id} \}$, the formula is given explicitly in [MM].
\par
\vskip 4pt
We deduce the Counting Formula 2.3 in \S 2.
\par
\vskip 4pt
We shall then give three applications.
\par
\vskip 4pt
First, as an immediate consequence of the Counting Formula,
we see that a K3 surface
$X$ such that $\rho(X) \geq 3$ and $det\, NS(X)$ is square free, has no
FM partner other than $X$ itself
(Corollary 2.7 3)). This is in sharp contrast to
the cases of $\rho = 1, 2$ described in [Og1]
(See also Corollary 2.7 4) and Theorem 3.3). We shall also show
how some of known cases
are derived from our Counting Formula (Corollary 2.7 1) and 4)).

Next, as an arithmetical application, we
recast two famous questions on class numbers, after Gauss, circa 1800, into geometric
terms of FM partners of K3 surfaces of Picard number 2 (Theorem 3.3,
Questions I' and II' in \S 3). Here we recall the two questions
on the class number $h(p)$ of the real
quadratic field $\mathbf Q(\sqrt{p})$ for a prime number
$p\equiv1~mod~4$:
\par
\vskip 4pt
\noindent
{\bf Question I.} {\it
Are there infinitely many primes $p$ such that $h(p)=1$?}
\par
\vskip 4pt
\noindent
{\bf Question II.} {\it
Is there a sequence of primes $p_1,p_2,...$ such that
$h(p_k)\rightarrow\infty$? }
\par
\vskip 4pt
Our reformulations (Questions I' and II') of these arithmetical questions
might afford
us a new geometrical approach to studying them via moduli
of stable sheaves on a family of K3 surfaces (cf. Theorem 2.2 and an
argument in \S 4).

Finally, for a more geometrical application, we shall construct
explicitly a pair
of smooth projective families
of K3 surfaces, $f : \mathcal X \rightarrow \Delta$ and
$g : \mathcal Y \rightarrow \Delta$, with a sequence
$\{t_{k}\}_{k=1}^{\infty} \subset \Delta$ such that
$\lim_{k \rightarrow \infty} t_{k} = 0$
and $\mathcal X_{t_{k}} \simeq \mathcal Y_{t_{k}}$, but
$\mathcal X_{0} \not\simeq \mathcal Y_{0}$ (Theorem 4.1). Here $\Delta$ is a
unit disk.
This explicitly shows the non-Hausdorff nature of
the moduli space of unpolarized K3 surfaces.
We carry out the construction in \S 4.

We remark that, aside from the three applications here,
the Counting Formula 2.3 has
also been applied to a study of Kummer structures on a K3 surface [HLOY2].
\par
\vskip 4pt
\noindent
{\bf Acknowledgement.}
We thank Professors B. Gross and K. Zuo for their interest in this work.
We also thank C. Madonna for communications on our previous note [HLOY1].
Main part of this collaboration was done during the first and third named
authors' stay at Harvard University. They would like to express their
thanks to the Harvard University and the Education Ministry of Japan
for financial support. The second named author is supported by
NSF grant DMS-0072158. The third named author would like to express his thanks 
to the organizers for their invitation to the conference.

\section{G-equivalence classes of primitive embeddings}

The goal of this section is Theorem 1.4. Here
we shall a need a $G$-equivariant version of
the notion of primitive embeddings [Ni1].

Unless otherwise stated, by a lattice $L := (L, (*, **))$,
we mean a pair of a free $\mathbf Z$-module $L$ of finite rank
and its non-degenerate symmetric integral-valued bilinear form
$(*, **)  : L \times L \rightarrow \mathbf Z$.
We set $L^*=\text{Hom}(L,\mathbf Z)$.

Let $\Lambda = (\Lambda, (*, **))$ be an even, unimodular, indefinite lattice.
Consider an even lattice $T$ which admits at least one primitive
embedding $\iota_{0}: T \hookrightarrow \Lambda$, and a subgroup $G
\subset O(T)$. The group $G$ is not necessarily a finite group. Let $S$ be a
lattice isomorphic to the orthogonal
lattice $\iota_{0}(T)^{\perp}$ in $\Lambda$
and choose an isomorphism $\iota_{0}^{S} : S \simeq \iota_{0}(T)^{\perp}$.
Throughout \S 1, the objects $\Lambda$,
$T$, $\iota_{0}$, $G$, $S$ and $\iota_{0}^{S}$ are fixed.

\begin{definition} \label{def:defTa}
Two primitive embeddings $\iota:
T\hookrightarrow \Lambda$, $\iota': T \hookrightarrow \Lambda$ are
called $G$-equivalent if there exist $\Phi \in O(\Lambda)$ and $g \in G$
such that the following diagram commutes:
$$
\begin{matrix}  \Lambda & \mapright{\Phi} & \Lambda &\cr
         \mapup{\iota} &  & \mapup{\iota'} &\cr
        T & \mapright{g} & T & \;\;.\cr
\end{matrix}
$$
\end{definition}

Our main interest in this section is the set of all primitive
embeddings:
$$
\PE:=\{ \iota: T \hookrightarrow \Lambda ; \text{primitive}\;
\text{embedding} \}\, ,
$$
and the set of $G$-equivalence classes, i.e. the quotient set
$$
\GPE :=\PE / \text{$G$-equivalence} \;.
$$
\par
We denote the restriction of the
bilinear form $(\;,\;)$ on $\Lambda$ to $S$ and $T$, respectively, by $(\;,\;)_S$,
$(\;,\;)_T$. These coincide with the original bilinear forms of $S$ and $T$.
(Here and hereafter, we often identify $T$ and $S$ with the primitive
sublattices $\iota_{0}(T)$ and $\iota_{0}(T)^{\perp} = \iota_{0}^{S}(S)$
of $\Lambda$ {\it via the fixed isomorphisms} $\iota_{0}$ and
$\iota_{0}^{S}$. )
Since $\Lambda$ is unimodular and $S$ and $T$ are primitive
in $\Lambda$, we have surjective maps $\pi_S: \Lambda = \Lambda^{*}
\rightarrow S^*$, $\pi_T: \Lambda = \Lambda^{*} \rightarrow T^*$
defined by $\pi_S(l)=(l,*)|_S$, $\pi_T(l)=(l,*)|_T$, and
natural isomorphisms induced by $\pi_{S}$ and $\pi_{T}$
([Ni1, Corollary 1.6.2]):
$$
\begin{matrix}
      &     \bar\pi_S  &                    &  \bar\pi_T     &        \cr
S^*/S & \simleftarrow  & \Lambda/(S\oplus T) & \simrightarrow & T^*/T
\quad . \cr
\end{matrix}
\label{eqn:eqIa}
$$
Here the last isomorphisms are written explicitly as
$\bar\pi_S(l \; \text{mod} \; S \oplus T  ) = l_S \; \text{mod} \; S$
and $\bar\pi_T(l \; \text{mod} \; S \oplus T  ) = l_T \; \text{mod} \; T$,
where $l_S,l_T$ are defined by the orthogonal decomposition $l=l_S+l_T$
made in $\Lambda \otimes \mathbf Q =(S \oplus T)\otimes \mathbf Q$. We make the identifications
$S \subset S^* \subset S\otimes\mathbf Q$  and
$T \subset T^* \subset T\otimes\mathbf Q$  via the respective nondegenerate
forms $(\;,\;)_S$ and $(\;,\;)_T$.

As it is defined by [Ni1, \S 1], the discriminant group
$(A_{S}, q_{S})$ of the lattice $(L,(\;,\;)_L)$ is the pair of the
abelian group $A_L:=L^*/L$ and a natural quadratic form
$q_L: A_L \rightarrow  \mathbf Q/ 2 \mathbf Z$ defined by
$
q_L(x \; \text{mod} \; L):=
(x, x)_{L \otimes \mathbf Q} \; \text{mod}\; 2\mathbf Z$.
Here $(\;,\;)_{L \otimes \mathbf Q}$ is the natural
$\mathbf Q$-linear extension
of the bilinear form $(\;,\;)_{L}$ of $L$.

\begin{proposition} \label{prop:propIa}{\rm ([Ni1, Corollary 1.6.2])}
Put $\sigma_0=\bar\pi_S\circ \bar\pi_T^{-1}:
A_T \simrightarrow A_S$.
Then $\sigma_0$ is an isometry:
$\sigma_0:(A_T,-q_{A_T}) \cong (A_S, q_{A_S})$.
\end{proposition}

As it is well-known, primitive embedding of $T$ into $\Lambda$ are
closely related to the so-called genus of the lattice $S$.
We shall briefly recall this here.

\begin{definitionThm}   \label{def:genusS}
Let $M$ be an even lattice.
The set of isomorphism classes of
the lattices $M'$ satisfying the following equivalent conditions
1) and 2), are called the {\it genus} of $M$,
and will be denoted by $\mathcal G(M)$:
\begin{list}{}{
\setlength{\leftmargin}{10pt}
\setlength{\labelwidth}{6pt}
}
\item[1)] $M'\otimes \mathbf Z_p \cong M \otimes \mathbf Z_p$ for all primes $p$,
          and $M' \otimes \mathbf R \cong M \otimes \mathbf R$.
\item[2)]  $(A_{M'}, q_{M'}) \cong (A_{M}, q_{M})$,
           and $\text{sgn}\,M' =\text{sgn}\,M$
where $\text{sgn}\, M$ is the signature of the lattice $M$.
\end{list}
\end{definitionThm}
The equivalence of 1) and 2) is due to [Ni1, Corollary 1.9.4]. From the
condition 2), we have
$\text{det}\, M' =\text{det}\, M$ and $\text{rank}\, M' =\text{rank}\, M$
if $M' \in \mathcal G (M)$. It is well-known that the genus $\mathcal G(M)$
is a finite set (see e.g. [Cs, Page 128,
Theorem 1.1]).

\vskip0.5cm

Consider the genus $\mathcal G(S)$ of $S$ and put:
$
\mathcal G(S)=\{ S_1,S_2,\cdots,S_m\} \;\;(S_1:=S)$.
By definition, we have an isomorphism $(A_{S},q_{S})\cong
(A_{S_j},q_{S_j})$ for each $S_j \in \mathcal G(S)$. Since
$\sigma_0: (A_T,-q_T) \cong (A_S,q_S)$ by Proposition {\ref{prop:propIa}},
we can then choose
an isomorphism for each $S_j \in \mathcal G(S)$:
$
\varphi_j: (A_T,-q_T) \simrightarrow (A_{S_j}, q_{S_j})$.

Throughout \S 1, these isomorphisms
$\varphi_{j}$ ($1 \leq j \leq m$) are fixed.
Now consider an arbitrarily primitive embedding
$\iota: T \hookrightarrow \Lambda$.
We have, by Proposition 1.2,
$$
(A_{\iota(T)^\perp},q_{\iota(T)^\perp}) \cong
(A_T,-q_T) \cong (A_S, q_S)
$$
and also
$\text{sgn}\, \iota(T)^\perp =\text{sgn}\, \Lambda -\text{sgn}\,T$.
Therefore there exists
a unique $j = j(\iota) \in \{1,  2,  \cdots , m\}$ and an isomorphism
$
S_j \simrightarrow \iota(T)^\perp$. If two primitive embeddings
$\iota:T \hookrightarrow \Lambda$ and
$\iota': T \hookrightarrow \Lambda$
are $G$-equivalent, there exist elements
$\Phi \in O(\Lambda)$ and $g \in G$ such that
$\Phi\circ \iota=\iota'\circ g$.
Since the restriction $\Phi|_{\iota(T)}$ is then an
isometry from $\iota(T)$ to $\iota'(T)$, we have
$S_{j(\iota)} \cong \iota(T)^\perp \cong \iota'(T)^\perp
\cong S_{j(\iota')}$, which
implies that $j(\iota) = j(\iota')$, i.e. $S_{j(\iota)} = S_{j(\iota')}$ in
$\mathcal G(S)$. Therefore setting
$$
\PEj :=\{ \iota: T \hookrightarrow \Lambda \in \PE \;| \;
\iota(T)^\perp \cong S_j \;\} ,
\label{eqn:pj}
$$
$$
\GPEj :=\{ [\iota: T \hookrightarrow \Lambda] \in \GPE \;| \;
\iota(T)^\perp \cong S_j \;\} ,
\label{eqn:pej}
$$
we have the following well-defined disjoint unions:
$$
\PE=\bigcup_{j=1}^m \PEj \;\;,\;\;
\GPE=\bigcup_{j=1}^m \GPEj \;\;.
$$
The following is a generalization of [MM], which will be proved in 
Appendix A:

\begin{theorem} \label{thm:numberTh}
For each $j$, we have
$\vert \GPEj \vert =  \vert O(S_j) \setminus O(A_{S_j})/G \vert$. In
particular,
$$
\vert \GPE | = \sum_{j=1}^m \vert O(S_j) \setminus O(A_{S_j})/G \vert\, .$$
Here $O(S_j)$ acts on $O(A_{S_j})$ through the natural map $O(S_j)
\rightarrow O(A_{S_j})$, $h \mapsto \bar h$. $G$ acts on $O(A_{S_j})$
through the composite of the natural map $G \subset O(T)
\rightarrow O(A_T)$
and the adjoint map $\text{ad}(\varphi_{j}) : O(A_{T})
\rightarrow O(A_{S_{j}})$.
\end{theorem}

\section{ FM number of a K3 surface}

Let $D(X)$ be the derived category of bounded complexes of coherent
sheaves on a smooth projective variety $X$.
A smooth projective variety $Y$ is called
a {\it Fourier-Mukai (FM) partner} of $X$ if there is an equivalence
$D(Y)\cong D(X)$ as triangulated categories.
We denote by $FM(X)$ the set of isomorphism classes of
FM partners of $X$.

\begin{definition} \label{def:dumb}
We call the cardinality $|FM(X)|$ of the set $FM(X)$
the Fourier-Mukai (FM) number of $X$.
\end{definition}

A FM partner of a K3 surface
is again a K3 surface by a result of Mukai [Mu3] and Orlov [Or]
(See also [BM]).

By a K3 surface, we mean a smooth {\it projective} surface $X$
over $\mathbf C$ with
$\mathcal O_X(K_X)\cong \mathcal O_X$ and $h^1(\mathcal O_X)=0$.
We denote by $(*,**)$ the
symmetric bilinear form on $H^2(X,\mathbf Z)$ given by the cup product. Then
$(H^2(X,\mathbf Z), (*,**))$ is an even unimodular lattice of signature
$(3,19)$. This lattice is isomorphic to the K3 lattice $\Lambda_{\text{K3}}:=
E_8(-1)^{\oplus 2}\oplus U^{\oplus 3}$ where $U$ is the hyperbolic lattice
(an even unimodular lattice of signature (1,1)).
An isomorphism $\tau: H^2(X,\mathbf Z)\rightarrow \Lambda_{\text{K3}}$
is called a marking and a pair
$(X,\tau)$ is called a marked K3 surface. We denote by
$NS(X)$ the N\'eron-Severi lattice of $X$ and by
$\rho(X)$ the Picard number, i.e. the rank of $NS(X)$. The lattice $NS(X)$
is primitive in $H^2(X,\mathbf Z)$ and has signature $(1,\rho(X)-1)$. We call
the orthogonal lattice $T(X):=NS(X)^\perp$ in $H^2(X,\mathbf Z)$ the
transcendental lattice. $T(X)$ is primitive
in $H^2(X,\mathbf Z)$ and has signature $(2,20-\rho(X))$. We denote by
$\omega_X$ a nowhere vanishing holomorphic two form on $X$. Then one
has a natural inclusion
$$
\mathbf C \omega_X \oplus \mathbf C \bar \omega_X \subset T(X)\otimes \mathbf C.
$$
This defines a Hodge structure $(T(X),\mathbf C \omega_X)$ of
weight 2 on $T(X)$.
Note that  $T(X)$ is the minimal primitive sublattice of
$H^{2}(X, \mathbf Z)$ such that
$T(X) \otimes \mathbf C$ contains $\mathbf C \omega_{X}$.

The following is a fundamental theorem due to Mukai and Orlov.
This connects three aspects of FM partners of K3 surfaces. Statement
1) is categorical, 2) is arithmetical, and 3) is geometrical.

\begin{theorem} \label{thm:Hodge} {\rm ([Mu3],[Or])}
Let $X$ and $Y$ be K3 surfaces. Then the following statements are equivalent:
\begin{list}{}{
\setlength{\leftmargin}{10pt}
\setlength{\labelwidth}{6pt}
}
\item[1)] $D(Y) \cong D(X)$;
\item[2)] There exists a Hodge isometry $g: (T(Y),\mathbf C \omega_Y)
\simrightarrow (T(X), \mathbf C \omega_X))$;
\item[3)] $Y$ is a two dimensional fine compact moduli space
of stable sheaves on $X$ with respect to some polarization on $X$.
\end{list}
\end{theorem}

{\it For the rest of this section, $X$ will be a K3 surface.} Our Counting Formula 2.3 follows from specializing Theorem 1.4 to the case
$$
\Lambda=H^2(X,\mathbf Z),~~~(T,\mathbf C\omega)=(T(X),\mathbf C\omega_X),
~~~S=NS(X),~~~G=O_{Hodge}(T(X),\mathbf C\omega_X).
$$
Note that in this case the group
$G$ is a finite cyclic group of order $2I$ such that
$\varphi(2I) \vert \text{rk}\, T(X)$. (See Appendix B for proof.) 

The aim of this section is to prove the following theorem:
\begin{theorem} \label{thm:FMcount} {\bf (The Counting Formula) }
Set $\mathcal G(NS(X))=\{ S_1=S,S_2,\cdots , S_m \}$
where $m=|\mathcal G(NS(X))|$. Then,
$$
|FM(X)| = \sum_{j=1}^m |O(S_j)\setminus O(A_{S_j}) /
O_{Hodge}(T(X),\mathbf C \omega_{X})|.
$$
Moreover the $j$-th summand here coincides with the number
of FM partners $Y\in FM(X)$ with $NS(Y)\simeq S_j$.
\end{theorem}

\noindent
{\it Proof.} By Theorem 1.4, it suffices to show the following:

\begin{theorem} \label{thm:FMtoPE}
There is a natural bijective correspondence
$$
\mathcal P^{G-eq}(T(X),H^2(X,\mathbf Z)) \leftrightarrow FM(X).
$$
Moreover, under this bijection a primitive embedding $\iota$ with
$\iota(T)^\perp\simeq S_j$ corresponds to a FM partner
$Y$ of $X$ with $NS(Y)\simeq S_j$.
\end{theorem}

Our proof of Theorem 2.4 is based on Theorem 2.2 and the
two most fundamental facts about K3 surfaces, namely, the global Torelli
theorem and the surjectivity of period map (See e.g. [BPV, Chapter VIII]
and the
references therein).

In what follows, for convenience, we identify
$H^2(X,\mathbf Z)=\Lambda_{K3}$ through a fixed marking
$\tau_X$ and set $(T,\mathbf C\omega)=(T(X),\mathbf C\omega_X)\subset H^2(X,\mathbf Z)$
under the natural inclusion.
We proceed in five steps.

\vskip .1in

\noindent
{\it Step 0.}
For an arbitrary primitive embedding $
\iota: T \hookrightarrow \LKIII$, there exists a marked K3 surface
$(Y,\tau_Y)$ with the commutative diagram:
$$
\begin{matrix} \tau_Y: & H^2(Y,\mathbf Z) & \simrightarrow & \LKIII \cr
                & \cup           &                & \mapup{\iota} \cr
\iota^{-1}\circ \tau_Y\vert_{T(Y)}:
          & (T(Y),\mathbf C \omega_Y) & \rightarrow & (T,\mathbf C \omega ) \cr
\end{matrix}
\label{eqn:FII}
$$
\begin{proof}
This follows immediately from the
surjectivity of the period map and the minimality of $T$.
\end{proof}


\noindent
{\it Step 1.} There is a map $
\tilde c:  \PEKIII  \rightarrow  FM(X)$, i.e. there exists a
FM partner $Y(\iota)$ for each primitive
embedding $\iota: T \hookrightarrow \LKIII$.

\begin{proof}
Let $\iota: T \hookrightarrow \LKIII$ be a primitive
embedding, then by Step 0. we have a marked K3 surface
$(Y,\tau_Y)$ with the commutative diagram in Step 0.
Then, by Theorem 2.2,
we have $Y \in FM(X)$. Assume there exists another marked K3
surface $(Y',\tau_{Y'})$
with the commutative diagram in Step 0 ($Y$ replaced by $Y'$). Then
$$
\tau_Y^{-1}\circ\tau_{Y'}:(H^2(Y',\mathbf Z), \mathbf C \omega_{Y'})
\rightarrow (H^2(Y,\mathbf Z) , \mathbf C \omega_Y)
$$
is a Hodge isometry, which means that $Y'\cong Y$ by the
Torelli theorem (See e.g. [BPV]). Therefore
there exists a well-defined map
$\tilde c:  \PEKIII \rightarrow  FM(X)$.
\end{proof}

\noindent
{\it Step 2.} The map $\tilde c$ is surjective.

\begin{proof} Let $Y \in
FM(X)$ and fix a marking $\tau_Y: H^2(Y,\mathbf Z) \rightarrow \LKIII$. Then
by Theorem 2.2,
we have a Hodge isometry $g: (T(Y),\mathbf C \omega_Y)
\simrightarrow (T(X),\mathbf C \omega_X)$ and the following
commutative diagram in which $T'=\tau_Y(T(Y)), \omega'=\tau_Y(\omega_Y)$:
$$
\begin{matrix}
\LKIII & \simmapleftarrow{\tau_Y} & H^2(Y,\mathbf Z) & 
& H^2(X,\mathbf Z) & \simmaprightarrow{\tau_X} & \LKIII \cr
\cup \; \iota_Y & & \cup &  & \cup &  & \cup\; \iota_X \cr
(T',\mathbf C \omega') & \simmapleftarrow{} & (T(Y), \mathbf C \omega_Y) &
\simmaprightarrow{g} & (T(X),\mathbf C \omega_X) & \simmaprightarrow{} &
 (T,\mathbf C \omega) \cr
\end{matrix}
$$
The map $
\iota:=\iota_Y\circ \tau_Y|_{T(Y)}\circ g^{-1}\circ\tau_X^{-1}|_{T}:
T \hookrightarrow \LKIII$
is then a primitive embedding and $Y = \tilde c(\iota)$.
\end{proof}

\noindent
{\it Step 3.} If two primitive embeddings $\iota,\iota':
(T,\mathbf C \omega) \hookrightarrow \LKIII$ are $G$-equivalent, then
$\tilde c(\iota) =\tilde c(\iota')$.
This together with the surjectivity
in Step 2 entails that the map $\tilde c$
descends to a surjective map $c: \GPEKIII \rightarrow FM(X)$.

\begin{proof}  Suppose $\iota,\iota':
(T,\mathbf C \omega) \hookrightarrow \LKIII$
are $G$-equivalent, and denote $Y=\tilde c(\iota), Y'=\tilde c(\iota')$.
We also fix markings $\tau_Y, \tau_{Y'}$ of $Y$ and $Y'$.
Then we have the following commutative diagram:
$$
\begin{matrix}
H^2(Y,\mathbf Z) & \simmaprightarrow{\tau_Y}  & \LKIII
               & \simmaprightarrow{\exists\Phi}    & \LKIII
               & \simmapleftarrow{\tau_{Y'}} & H^2(Y',\mathbf Z) \cr
\cup           &                    & \uparrow{\iota}
               &                    & \uparrow{\iota'}
               &                    & \cup \cr
(T(Y),\mathbf C \omega_Y) & \simmaprightarrow{}  & (T, \mathbf C \omega)
               & \simmaprightarrow{\exists g} & (T,\mathbf C \omega)
               & \simmapleftarrow{}         & (T(Y'), \mathbf C \omega_{Y'}) \cr
\end{matrix}
$$
where $\Phi\in O(\LKIII)$ and $g \in G=O_{Hodge}(T, \mathbf C \omega)$.
{}From this diagram, we have a Hodge isometry $\tau_{Y'}^{-1}\circ
\Phi \circ \tau_Y: (H^2(Y,\mathbf Z),\mathbf C \omega_Y)
\rightarrow (H^2(Y',\mathbf Z),\mathbf C \omega_{Y'})$, which implies
$Y\cong Y'$ by the global Torelli theorem.
\end{proof}

\noindent
{\it Step 4.} The induced map $c: \GPEKIII \rightarrow FM(X)$
is injective.

\begin{proof} Take two $G$-equivalence classes $[\iota],
[\iota'] \in \GPEKIII$ and denote
$Y=c([\iota])=\tilde c(\iota)$, $Y'=c([\iota'])=\tilde c(\iota')$. We fix
markings $\tau_Y, \tau_{Y'}$ for $Y$ and $Y'$, respectively. To see the
injectivity we show that $f:Y'\cong Y$ implies $[\iota']=[\iota]$. Denote
the Hodge isometry $f^*: (H^2(Y,\mathbf Z),\mathbf C \omega_Y)
\rightarrow (H^2(Y',\mathbf Z),\mathbf C \omega_{Y'})$ associated with the
isomorphism $f: Y' \cong Y$. Then we have the following commutative diagram:
$$
\begin{matrix}
\LKIII & \simmapleftarrow{\tau_Y} & H^2(Y,\mathbf Z) &
       \simmaprightarrow{f^*}
       & H^2(Y',\mathbf Z)  & \simmaprightarrow{\tau_{Y'}} & \LKIII \cr
\mapupleft{\iota \; } &  & \cup &  & \cup &  & \mapup{\iota'} \cr
(T,\mathbf C \omega) & \mapleft{\sim} & (T(Y), \mathbf C \omega_Y) &
 \mapright{\sim} & (T(Y'),\mathbf \omega_{Y'}) &
                \mapright{\sim} & (T,\mathbf C \omega) \cr
\end{matrix}
$$
{}From this diagram we obtain;
$$
\begin{matrix}
\LKIII & \simmaprightarrow{\exists\Phi} & \LKIII \cr
\mapupleft{ \iota \; } &  & \mapup{\iota'} \cr
(T,\mathbf C\omega) & \simmaprightarrow{\exists g} & (T,\mathbf C \omega) \cr
\end{matrix}
$$
where $g \in O_{Hodge}(T,\mathbf C \omega)$. We then conclude
$[\iota]=[\iota']$ in $\GPEKIII$.
\end{proof}

It now follows that $c$ is bijective.
The second assertion in the theorem is clear from the construction of $c$.
\hfill $\square$

\vskip0.3cm

In an actual application of the Counting Formula 2.3,
we need to determine the genus of $NS(X)$. The following theorem due to
[Ni1, Theorem 1.14.2] is very useful when $\rho(X)$ is large:

\begin{theorem} \label{thm:nikulin}
Let $S$ be an even indefinite non-degenerate
lattice.
If $rk\, S \geq 2 + l(S)$, then $\mathcal G(S) = \{S\}$ and the natural map
$O(S) \rightarrow O(A_{S})$ is surjective. Here $l(S)$ is the minimal number
of generators of $A_{S}$, e.g. $l(A_{S}) = 1$ when $A_{S}$ is cyclic.
\end{theorem}

\begin{corollary} \label{prop:uniqueI}
If $\rho(X) \geq l(NS(X)) + 2$, then $FM(X) = \{X\}$.
\end{corollary}

\begin{proof}
By Theorem 2.5, we have $\mathcal G(NS(X)) = \{NS(X)\}$. Then, by the Counting
Formula, we have
$
|FM(X)| = |O(NS(X))\setminus O(A_{NS(X)}) /
O_{Hodge}(T(X),\mathbf C \omega_{X})|$.
Since the natural map $O(NS(X)) \rightarrow O(A_{NS(X)})$ is surjective
again by Theorem 2.5, we have $\vert FM(X) \vert = 1$.
\end{proof}

\begin{corollary} \label{cor:uniqueII}
\begin{list}{}{
\setlength{\leftmargin}{10pt}
\setlength{\labelwidth}{6pt}
}
\item[1)] {\rm (cf. [Mu1, Proposition 6.2])} If $\rho(X) \geq 12$,
then $FM(X) = \{X\}$. In particular, $FM(Km\, A) = \{Km\, A\}$.
\item[2)] If $\rho(X) \geq 3$ and $det NS(X)$ is square-free,
then $FM(X) = \{X\}$.
\item[3)] If $X \rightarrow \mathbf P^{1}$ is a Jacobian K3 surface, i.e.
an elliptic K3 surface with a section, then $FM(X) = \{X\}$.
\item[4)] {\rm (cf. [Og1, Proposition 1.10])} Assume that $\rho(X) = 1$ and
$NS(X) = \mathbf Z H$ with
$(H^{2}) = 2n$. Then $\vert FM(X) \vert = 2^{\tau(n) -1}$. Here $\tau(1)=1$,
and
$\tau(n)$ is the number of prime factors of $n\geq2$, e.g.
$\tau(4) = \tau(2) = 1$, $\tau(6) = 2$.
\end{list}
\end{corollary}

\noindent{\it Proof.}
\par
\vskip 4pt
\noindent
{\it 1):} This follows from Corollary 2.6 and the following calculation:
$l(NS(X))  = l(T_{X}) \leq rk\, T(X) \leq 10 \leq\, \rho(X) - 2$.
\par
\vskip 4pt
\noindent
{\it 2):} Since $\text{det}\, NS(X)$ is square-free, it follows that
$A_{NS(X)}$ is cyclic. Then $l(NS(X)) = 1$. Since $\rho(X) \ge 3$, the result
follows from Corollary 2.6.
\par
\vskip 4pt
\noindent
{\it 3):} The unimodular hyperbolic lattice $U$ is a direct summand
of $NS(X)$, say $NS(X) = U \oplus N$, by the assumption.
Thus, $l(NS(X)) = l(N) \leq rk\, N = \,
\rho(X) - 2$. Therefore, the result follows from Corollary 2.6.
\par
\vskip 4pt
\noindent
{\it 4):} We assume $n \geq 2$. The case $n =1$ is similar
and easier.
Since $NS(X) = \mathbf Z H$ with $(H^{2}) = 2n$, we have
$\mathcal G (NS(X)) = \{NS(X)\}$, $O(NS(X)) \simeq O_{Hodge}(T(X), \mathbf C \omega_{X})) \simeq \mathbf Z /2$ and
$(A_{NS(X)}, q_{NS(X)}) \simeq (\mathbf Z /2n, q_{2n})$. Here $q_{N}$ on
$\mathbf Z /N$ is defined by $q_{N}(1) = 1/N$.
Therefore,
$\vert FM(X) \vert = \vert O(\mathbf Z /2n, q_{2n}) \vert /2$
by the Counting Formula. By a straightforward calculation (cf. [Og1]),
we have also $\vert O(\mathbf Z/2n, q_{2n}) \vert = 2^{\tau(n)}$. \qed

\section{ Class Numbers and Fourier-Mukai Numbers}

In this section, we focus on K3 surfaces of $\rho(X) = 2$ (Theorem 3.3) and
recast the two class number problems, Questions I and II
in Introduction, in corresponding geometrical terms (Questions I' and II').

We begin by recalling the following fact which is well-known to experts
(see [Ni 1, 2] and [Mo]).
This proposition shows that every question about
the class numbers of real quadratic fields can, in principle,
be recast in terms of K3 surfaces.

\begin{proposition} \label{prop:ehl}
For any even hyperbolic lattice with $rk~S\leq 10$, there is a projective
K3 surface $X$ with $NS(X)\simeq S$. In particular,
for every rank 2 even hyperbolic lattice $S$, there is a projective K3
surface $X$ with $NS(X)\simeq S$.
\end{proposition}

Next we briefly recall some basic background on class numbers from
[Cs] and [Za]. Let $D$ be an odd fundamental discriminant, i.e. a rational
integer
$D$ such that $D\equiv1~mod~4$ and $D$ square-free. Put $K:=\mathbf Q(\sqrt D)$, and let $\mathcal O_K$
denote its ring of integers. Throughout this section, this $D$ is fixed.
Let $H(D)$ denote the narrow ideal class group of $\mathcal O_K$, i.e.
the quotient of the group of fractional ideals of $\mathcal O_K$,
by the subgroup of principal fractional ideals with {\it positive} norm.
A theorem of Hurwitz says that $H(D)$ is finite.
The number $h(D)=\vert H(D)\vert$ is {\it the class number} of $K$.

\begin{theorem}  \label{thm:Bijection} {\rm ([Za], [Cs])}
The following three sets are in natural 1-1 correspondence:
\begin{list}{}{
\setlength{\leftmargin}{10pt}
\setlength{\labelwidth}{6pt}
}
\item[1)] The class group $H(D)$.
\item[2)] The set $\mathcal B(D)$ of proper, i.e.
$SL(2,\mathbf Z)$-equivalence classes of integral binary quadratic forms
$a x^2+bxy+cy^2$ with discriminant $b^2-4ac=D$.
\item[3)] The set $\mathcal L(D)$ of proper, i.e. orientation-preserving,
isomorphism classes
of even hyperbolic rank 2 lattices $S$ with $det~S=-D$.
\end{list}
\end{theorem}

Note that every class in $\mathcal L(D)$ can be represented by a pair
$\left(\mathbf Z^2,
\left(
\begin{matrix}
2a&b\cr b&2c
\end{matrix}
\right)\right)$, where the matrix represents
the bilinear form with respect to the standard base of $\mathbf Z^2$.
Two such pairs are properly isomorphic if and only if
their matrices are conjugate under $SL(2,\mathbf Z)$.

We recall a few more facts about the sets
$H(D),\mathcal B(D), \mathcal L(D)$.
For details, see [Za, \S12] and
[Cs, Chap. 14 \S3].
Consider the exact sequence
$$
0\rightarrow\mathcal J\rightarrow H(D)\mapright{sq}H(D)^{2} \rightarrow 0
$$
where $sq(x) = x^{2}$ and $\mathcal J$ is the kernel of $sq$.
It is known that
$\vert\mathcal J\vert = \vert H(D)/H(D)^2\vert = 2^{n-1}$,
where $n$ is the number of prime factors of $D$.
Given a rank 2 even lattice $L$ with $det~L=-D$, let
$\tilde{\mathcal G}(L)$ be the set of {\it proper}
isomorphism classes of lattices in the same genus as $L$.
The set of genera $\{\tilde{\mathcal G}(L)| [L]\in\mathcal L(D)\}$
is a principal homegeneous space of the group $H(D)/H(D)^2$.
Therefore, there are $[L_1],...,[L_{2^{n-1}}]\in\mathcal L(D)$,
and a decomposition:
$$
\mathcal L(D)=\coprod_{k=1}^{2^{n-1}}\tilde{\mathcal G}(L_k)
$$
where the genera $\tilde{\mathcal G}(L_k)$ all have the same cardinality.
Now the subset $\mathcal J\subset H(D)$ consists of the ideal classes
of the so-called 2-sided ideals. This means that,
under the correspondence $H(D)\leftrightarrow\mathcal L(D)$,
each element in $\mathcal J$ corresponds to a class in $\mathcal L(D)$
representable as a pair
$\left(\mathbf Z^2,
\left(
\begin{matrix}
2a&b\cr b&2c
\end{matrix}
\right)\right)$ such that
$
\left(\mathbf Z^2,
\left(
\begin{matrix}
2a&b\cr b&2c
\end{matrix}
\right)\right)
\simeq
\left(\mathbf Z^2,
\left(
\begin{matrix}
2a&-b\cr -b&2c
\end{matrix}
\right)\right)
$
under $SL(2,\mathbf Z)$.
In particular the number of such classes in $\mathcal L(D)$ is
$\vert\mathcal J\vert=2^{n-1}$.

Now we have:
\begin{theorem} \label{thm:MainII}
If $D=p$ is a prime and $X$ is a K3 surface with
$\rho(X)=2$ and $det~NS(X)=-p$,
then
$\vert FM(X)\vert= (h(p)+1)/2$.
\end{theorem}

\vskip0.3cm
\noindent
{\bf Remark.}
For each prime number $p\equiv1~mod~4$, there is a K3 surface
$X$ such that $NS(X)\simeq S_0$ by Proposition 3.1. Here
$S_0=\left(\mathbf Z^2,
\left(
\begin{matrix}
2&1\cr 1&\frac{1-p}{2}
\end{matrix}
\right)\right)$. Note that $S_{0}$ is an even hyperbolic lattice of rank $2$
and of $det\, S_{0} = -p$.
$\qed$

\begin{proof}
Put $S=NS(X)$. Since $D=p$, i.e. $n=1$, it follows that
$\tilde{\mathcal G}(S)=\mathcal L(D)$, hence
$
h(p)=\vert \tilde{\mathcal G}(S)\vert$. Since $\mathcal G(S)$ consists of
isomorphism (not necessarily proper)
classes of lattices, it follows that $\mathcal G(S)$ is the
orbit space of the group
$\langle\sigma\rangle=GL(2,\mathbf Z)/SL(2,\mathbf Z) (\simeq C_{2})$
acting on $\tilde{\mathcal G}(S)$.
Since $\vert\mathcal J\vert=1$ because $n=1$,
it follows that there is just one class $[L]\in\tilde{\mathcal G}(S)$
(in fact, $[L]=[S_0]$) which is fixed by $\sigma$
acting on $\tilde{\mathcal G}(S)$. Since $\sigma$ is an involution, it follows that $
\vert\mathcal G(S)\vert= (\vert\tilde{\mathcal G}(S)\vert+1)/2$.
On the other hand,
the group $A_S$ has $p$ elements.
Since $p$ is prime, it follows that $A_S\simeq\mathbf Z/p$,
and that $O(A_S)=\{\pm id\}$. Since the group $O_{Hodge}(T(X),\mathbf Z\omega_X)$
obviously contains $\pm id$, it follows that each of the quotients
$O(A_{S_i})/O_{Hodge}(T(X),\mathbf C\omega_X)$ in
the Counting Formula for $\vert FM(X)\vert$
has just 1 element. Therefore
$\vert FM(X)\vert=\vert\mathcal G(S)\vert.$
Now we are done.
\end{proof}

{}From a table [Fl], for a prime $1297$ ($\equiv 1\, \text{mod}\, 4$), we have
$h(1297) = 11$. Then for any K3 surface $X$ with $\rho(X) = 2$ and
$det\, NS(X) = -1297$, we have $\vert FM(X) \vert = 6$. In this case the value
$\vert FM(X) \vert$ is not a power of $2$, either. This example may be contrasted to the case $\rho(X) = 1$, where $\vert FM(X) \vert$ is always a power of
$2$
(Corollary 2.7 4).)
\par

\vskip 4pt

$$
\def\vspace#1{  height #1 & \omit && \omit && \omit &  \cr }
\vbox{\offinterlineskip
\hrule
\halign{ \strut
\vrule#&  $\;$ \hfil $#$ \hfil
&\vrule#&  $\;$ \hfil $#$ \hfil
&\vrule#&  $\;$ \hfil $#$ \hfil
&\vrule#
\cr
&  p &&  h(p)  &&  \vert FM(X)\vert  &\cr
\noalign{\hrule}
\noalign{\hrule}
& 229
&& 3
&& 2
&\cr
& 257
&& 3
&& 2
&\cr
& 401
&& 5
&& 3
&\cr
& 577
&& 7
&& 4
&\cr
& 733
&& 3
&& 2
&\cr
& 761
&& 3
&& 2
&\cr
& 1009
&& 7
&& 4
&\cr
& 1093
&& 5
&& 3
&\cr
& 1129
&& 9
&& 5
&\cr
& 1229
&& 3
&& 2
&\cr
& 1297
&& 11
&& 6
&\cr
& 1373
&& 3
&& 2
&\cr
& 1429
&& 5
&& 3
&\cr
& 1489
&& 3
&& 2
&\cr
}
\hrule}
$$

By Theorem 3.3, Gauss' questions can now be stated
in the following geometrical form:

\noindent {\bf Question I'.}\, {\it
Are there infinitely many
K3 surfaces $X$, with each $X$ having $\rho(X)=2$ and $FM(X)=\{X\}$,
such that the numbers $|det~NS(X)|$ are distinct primes?  }

\noindent {\bf Question II'.}\,
{\it
Is there a sequence of K3 surfaces $X_1,X_2,...$ with
$\rho(X_k)=2$ and $\vert FM(X_k)\vert\rightarrow\infty$
such that the numbers $|det~NS(X_k)|$ are primes?  }

\section{Non-Hausdorff properties of the moduli of K3 surfaces}
In this section, we shall show the following:
\begin{theorem}  \label{thm:nHdf}
Let $\Delta$ be a unit disk in $\mathbf C$. Then there
is a pair of smooth projective families of K3 surfaces
$\mathcal X \rightarrow \Delta$, $\mathcal Y \rightarrow \Delta$
and a sequence $\{t_{k}\} \subset \Delta$ such that
$$\lim_{k \rightarrow \infty} t_{k} = 0\, ,\, \mathcal X_{t_{k}}
\simeq \mathcal Y_{t_{k}}\, {\text but}\,
\mathcal X_{0} \not\simeq \mathcal Y_{0}\, .$$
\end{theorem}

In our construction, we shall make use of the 3-rd characterization of FM partners of a K3 surface
in Theorem 2.2. We begin by recalling a few facts about the $3$-rd characterization,
i.e. $2$-dimensional fine compact moduli spaces of stable sheaves on a K3 surface $X$.
\par
\vskip 4pt
Set $\widetilde{NS}(X):=H^0(X,\mathbf Z)
\oplus NS(X) \oplus H^4(X,\mathbf Z)$. We call $\widetilde{NS}(X)$
the extended N\'eron-Severi lattice of $X$. For
$v=(r,H,s) \in \widetilde{NS}(X)$ and an ample class $A \in NS(X)$,
we denote by $M_A(v)$ (resp. $\overline{M}_A(v)$)
the coarse moduli space of stable sheaves (resp. the coarse moduli space
of $S$-equivalence classes of
semi-stable sheaves) $\mathcal F$ with respect to the polarization $A$ with $\mu(\mathcal F)=v$.
Here $\mu(\mathcal F)$ is the so-called Mukai vector of $\mathcal F$,
which is defined by
$\mu(\mathcal F):=\text{ch}(\mathcal F) \sqrt{td_X}$.
By [Ma1, 2] (see also [Mu3]),  the space $\overline{M}_A(v)$ is a projective
compactification of $M_A(v)$, and
$\overline{M}_A(v) = M_A(v)$ if all the semi-stable sheaf $\mathcal F$
with $\mu(\mathcal F) = v$ is stable, for instance if $(r, s) = 1$.
\par
\vskip 4pt
The following theorem essentially due to Mukai [Mu3] is proved in [Or]
in the course of the proof of Theorem 2.2:

\begin{theorem}  \label{thm:finem}
If $Y \in FM(X)$, then $Y \simeq M_{H}((r, H, s))$ where
$H$ is ample, $r > 0$ and $s$ are integers
such that $(r, s) = 1$ and $2rs = (H^{2})$, and vice versa.
\end{theorem}

In our proof, we need is the following relative version:

\begin{lemma}   \label{lemma:famFM}
Let $\pi : (\mathcal X, \mathcal H) \rightarrow \mathcal B$ be a smooth
projective family of K3 surfaces.
Let $f : \overline{\mathcal M} \rightarrow \mathcal B$ be a relative moduli
space of the $S$-equivalence classes of semi-stable sheaves
of $\pi$ with respect to the polarization $\mathcal H$ with Mukai
vectors $(r, \mathcal H_{t}, s)$.
If $r > 0$, $(r, s) = 1$ and $2rs = (\mathcal H_{t}^{2})$, then $f$ is
projective and gives a relative FM family of $\pi$, i.e.
$\overline{\mathcal M}_{t} \in FM(\mathcal X_{t})$ for all
$t \in \mathcal B$.
\end{lemma}

\begin{proof} The existence of $f$ and its projectivity over
$\mathcal B$ are shown by Maruyama [Ma1] and [Ma2].
By the definition, we have
$\overline{\mathcal M}_{t} = \overline{\mathcal M}_{\mathcal H_{t}}((r,\mathcal H_{t},s))$.
Since $(r, s) = 1$, $2rs = (\mathcal H_{t}^{2})$ and $\mathcal H_{t}$ is ample on $\mathcal X_{t}$,
it follows that $\overline{\mathcal M}_{t} = \mathcal M_{H_{t}}((r,\mathcal H_{t},s)) \in FM(\mathcal X_{t})$ by Theorem 4.2.
\end{proof}
\noindent
{\it Proof of Theorem 4.1.}
Consider the number
$$m := \text{max}\, \{\rho(X) \vert X\,
\text{is a K3 surface with}\, \vert FM(X) \vert \geq 2\}\, .$$
Such a number $m$ exists and satisfies $m \leq 11$ by Corollary 2.7 1) and 4).
Take one such $X$ and choose then $Y \in FM(X)$ such that
$Y \not\simeq X$.
Write $Y = M_{H}(r, H, s)$ as in Theorem 4.2.
\par
\vskip 4pt
Let us consider the Kuranishi family
$u : \mathcal U \rightarrow \mathcal K$ of $X$. The base space
$\mathcal K$ is assumed to be a small polydisk in $H^{1}(X, T_{X}) \simeq
\bold C^{20}$. Choosing a marking $\tau
: R^{2}u_{*} \bold Z_{\mathcal U} \simeq \Lambda \times \mathcal K$, where
$\Lambda$ is the K3 lattice, we consider the period map
$$\pi : \mathcal K \rightarrow \mathcal P = \{ [\omega] \in \bold P(\Lambda\otimes\bold C)
\vert (\omega.\omega) = 0, (\omega.\overline{\omega}) > 0 \}
\subset  \bold P(\Lambda \otimes \bold C) \simeq \bold P^{21}\, .$$
By the local Torelli theorem (see e.g. [BPV, Chapter VIII]), $\mathcal K$ is
isomorphic, via $\pi$, to an open neighborhood in $\mathcal P$,
which we again denote by $\mathcal P$. In what follows, we identify
$\mathcal K$ with $\mathcal P$.
\par
\vskip 4pt
Set $\Lambda_{0} := \tau_{0}(NS(X))$.
Let us consider the subset $\mathcal K^{0}$ of $\mathcal K = \mathcal P$
defined
by the linear equations $([\omega], \Lambda_{0}) = 0$.
This $\mathcal K^{0}$ parametrizes K3 surfaces $\mathcal U_{t}$
such that $\Lambda_{0} \subset \tau_{t}(NS(\mathcal U_{t}))$,
or in other words, $\mathcal K^{0}$ is the subspace of $\mathcal K$ in which
$NS(X)$ is kept to be
algebraic.
Since $rank~\Lambda_0=\rho(X) = m$, we have $\text{dim}\, \mathcal K^{0} = 20 -m \geq 1$.
Let us take a generic small one-dimensional disk $\Delta$
such that $0 \in \Delta \subset \mathcal K^{0}$
and consider the family
$\varphi : \mathcal X \rightarrow \Delta$
induced from $u : \mathcal U \rightarrow \mathcal K$.
By definition, there is an invertible sheaf $\mathcal H$ on $\mathcal X$
such that $\mathcal H_{0} = H$. Since ampleness is an open condition,
$\mathcal H$ is $\varphi$-ample.
Therefore, the morphism $\varphi$ projective.
By [Og2, Main Theorem], we see that $NS(\mathcal X_{t}) \simeq \Lambda_{0} \simeq NS(X)$
for generic $t \in \Delta$
and that the set
$\mathcal S := \{s \in \Delta \vert \rho(\mathcal X_{s}) > m\}$
is everywhere dense in $\Delta$ with respect to the Euclidean topology.
\par
\vskip 4pt
Let us consider the family of FM partners $f : \mathcal Y \rightarrow \Delta$
of
$\varphi : \mathcal X \rightarrow \Delta$ with Mukai vectors
$(r, \mathcal H_{t}, s)$.
This exists by Lemma 4.3.
By the definition of $m$, $\mathcal S$ and $f$, one has
$\mathcal Y_{0} \simeq Y \not\simeq X \simeq \mathcal X_{0}$ and
$\mathcal Y_{s} \simeq \mathcal X_{s}$ for $s \in \mathcal S$.
Then, any sequence $\{t_{k}\} \subset \mathcal S$ with
$\lim_{k \rightarrow \infty} t_{k} = 0$ satisfies
the desired property.  \qed

\appendix
\section{Proof of Theorem 1.4}

In this appendix, for completeness, we shall give a sketch
of a proof of Theorem 1.4. Throughout this appendix, we use the same notation
as in \S 1.

Let $T$ and $S$ be the even lattices defined at the beginning of Section 1 and
fix $S_j \in \mathcal G(S)$. The goal is to show the
following equality:
$$
\vert \GPEj \vert =  \vert O(S_j) \setminus O(A_{S_j})/G \vert\, ,
$$
which clearly implies the first equality in Theorem {\ref{thm:numberTh}}.
This equality is a consequence of Propostions A.7 and A.9 below.

\begin{definition} \label{def:defIIa} {\rm ([Ni1, \S 1])}
An over-lattice
$L = (L,  (*, **)_{L})$ of $S_j \oplus T$ is
a pair consisting of a
$\mathbf Z$-module such that
$
S_j \oplus T \subset L \subset S_j^* \oplus T^*
$
and the bilinear form $(\;,\;)_L$ given by the restriction of
$(\;,\;)_{S_j^*}\oplus (\;,\;)_{T^*}$ to $L$. An over-lattice $L$ is
called {\it integral} if
its bilinear form $(\;,\;)_L$ is $\mathbf Z$-valued.
An integral over-lattice is a lattice in our sense.
\end{definition}

\begin{definition} \label{def:defIIb}
We denote by $\mathcal U_j:=\mathcal U(S_j\oplus T)$
the set of over-lattices $L = (L,(\;,\;)_L)$ of $S_j\oplus T$ such that
$L$ is integral, unimodular and such that the inclusion $T \subset L$ is
primitive. We also denote by $\mathcal E_j$ the subset of
$\mathcal U_j$ which consists of even lattices.
\end{definition}

\begin{lemma} \label{lemma:lemmaIIa} {\rm ([Ni1, \S 1])}
Let $L \in \mathcal U_j$. Then:
\begin{list}{}{
\setlength{\leftmargin}{10pt}
\setlength{\labelwidth}{6pt}
}
\item[1)] The orthogonal lattice $T^\perp$ in $L$ coincides with $S_j$, and
the inclusion $S_j \subset L$ is primitive.
\item[2)] The natural projections $\bar \pi_{L,S_j}$ and
$\bar \pi_{L,T}$ (see \S 1 for definition) are group isomorphisms.
\item[3)] $L$ is even, i.e. $L\in \mathcal E_j$ if and only if
$q_{S_j}(\varphi_L(x))=-q_T(x) \quad (\forall x \in A_T)$,
where $\varphi_L:=\bar\pi_{L,S_j}\circ \bar\pi_{L,T}^{-1}:
A_T \simrightarrow A_{S_j}$.
\end{list}
\end{lemma}

\begin{proof}
This follows directly from Proposition 1.2.
See also [Ni1, \S 1].
\end{proof}

\noindent
{\bf Remark.} 1) The relation
$q_{S_j}(\varphi_L(x))=-q_T(x) \quad (\forall x \in A_T)$ is nothing but
that the
isomorphism $\varphi_L$ is in fact an isometry of the discriminant lattices
$$
\varphi_L:(A_T,-q_T) \simrightarrow (A_{S_j},q_{S_j}) \;.
$$
\noindent
2) The isometry $\varphi_L:(A_T,-q_T) \simrightarrow (A_{S_j}, q_{S_j})$
associated to $L \in \mathcal E_j$ recovers $L$ explicitly as:
$
L/(S_j\oplus T)=\{ \varphi_L(a)\oplus a \vert a \in A_T \}
\subset S_j^*/S_j \oplus T^*/T$.   \qed

\begin{lemma} \label{lemma:LEj}
For an arbitrary isometry $f: (A_T,-q_T)
\simrightarrow (A_{S_j},q_{S_j})$, we define the over-lattice
$L^f$ of $S_{j} \oplus T$ by:
$
L^f/(S_j\oplus T)=\{ f(a) \oplus a \vert a \in A_T \}$.
Then $L^f$ is an element of $\mathcal E_j$, and the associated isometry
$\varphi_{L^f}(=\bar\pi_{L^f,S_j}\circ\bar\pi_{L^f,T}^{-1})$ coincides with
$f$.
\end{lemma}

\begin{proof} This is an old result known as gluing methods, see
[CoS, Chap. 4]
for example, and
also follows immediately from Proposition 1.2.
\end{proof}

By Lemma A.3, 3), each $L \in \mathcal E_j$ gives rise to an isometry
$
\varphi_L:(A_T,-q_T)\simrightarrow (A_{S_j},q_{S_j})$.
Combininig this with Proposition 1.2,
we obtain the following bijective correspondence
such that $L \mapsto \varphi_L, f \mapsto L^f$ are inverses of each
other:
$$
\mathcal E_j=\mathcal E(S_j\oplus T)  \leftrightarrow
\{ \varphi: (A_T,-q_T) \simrightarrow (A_{S_j},q_{S_j}) \}\, .
\label{eqn:corres}
$$

\begin{lemma} \label{lemma:OneToOne}
Fix an isometry $\varphi_j: (A_T,-q_T)
\simrightarrow (A_{S_j},q_{S_j})$. Then the following map is bijective:
$
\mu_j: O(A_{S_j}) \rightarrow \mathcal E_j\, ;\, \sigma \mapsto
L^{\sigma \circ \varphi_{j}}$.
\end{lemma}
\begin{proof}
The assignment $\sigma \mapsto \sigma \circ \varphi_{j}$ gives a
bijection between $O(A_{S_{j}})$ and the set
$\{ \varphi: (A_T,-q_T) \simrightarrow (A_{S_j},q_{S_j}) \}$.
Combining this bijection with Lemma A.3, we obtain the result.
\end{proof}

\begin{definition} \label{def:GequivL}
1) Two over-lattices $L, L' \in \mathcal E_j$
are called $G$-equivalent if there exist isometries $\Phi: L \simrightarrow
L'$ and $g: T \simrightarrow T$ such that $g \in G$ and $\Phi|_T=g$:
$$
\begin{matrix}  \Phi: & L & \simrightarrow &  L' \cr
               & \cup &  & \cup \cr
         g:    & T  & \simrightarrow & T \cr
\end{matrix}
$$
2) We denote the set of $G$-equivalence classes of the elements of
$\mathcal E_j$ by $\GEj$.
\end{definition}

\noindent
{\bf Remark.} By Lemma A.3, 1), the orthogonal lattice
$T^\perp$ in
$L$ coincides with $S_j$ for all $L \in \mathcal U_j$. Therefore when $L$
and $L'$
are $G$-equivalent, the isometry
$\Phi: L \simrightarrow L'$ determines $f:=\Phi_{S_j} \in O(S_j)$ as
well as $g=\Phi|_T \in G$. Conversely any pair $(f,g) \in O(S_j) \times G$
induces an isometry
$
f\oplus g : S_j^*\oplus T^* \simrightarrow S_j^*\oplus T^*$.
Moreover,  this isometry sends an over-lattice $L \subset S^*_j \oplus T^*$
to an
over-lattice $L'=(f\oplus g)(L) \subset S^*_j \oplus T^*$. Clearly $L$
and $L'$ are $G$-equivalent.  \qed

\begin{proposition} \label{prop:dcoset}
The bijective map $\mu_j: O(A_{S_j})
\rightarrow \mathcal E_j$ in Lemma {\ref{lemma:OneToOne}}
descends to the bijective map:
$
\bar\mu_j: O(S_j)\setminus O(A_{S_j}) / G \rightarrow
\GEj$,
where the double coset is defined as in Theorem {\ref{thm:numberTh}}.
\end{proposition}

\begin{proof}
Let $\sigma_1, \sigma_2 \in O(A_{S_j})$.  Assume
$L^{\sigma_1\circ \varphi_j}$ and $L^{\sigma_2\circ\varphi_j}$ are
$G$-equivalent over-lattices of $S_j\oplus T$, where $\varphi_j:(A_T,-q_T)
\simrightarrow (A_{S_j},q_{S_j})$ is a fixed isometry. By definition
of $G$-equivalence, we have isometries $\Phi: L^{\sigma_1\circ\varphi_j}
\simrightarrow L^{\sigma_2\circ\varphi_j}$ and $g: T \simrightarrow T$ such
that $\Phi|_T=g \in G$. By Lemma A.3, 1), $S_j=T^\perp$
in both $L^{\sigma_1\circ\varphi_j}$ and $L^{\sigma_2\circ\varphi_j}$,
and thus we have $\Phi(S_j\oplus T)=S_j\oplus T$ and $\Phi_{S_j}:=
\Phi\vert_{S_j} \in O(S_j)$, $\Phi_T:=\Phi\vert_T =g \in G$. Now we have
$
L^{\sigma_2\circ\varphi_j}/(S_j\oplus T) =
\{ \sigma_2\circ\varphi_j(a)\oplus a \; \vert \; a \in A_T \;\}$,
and
$$
\begin{aligned}
\Phi(L^{\sigma_1\circ\varphi_j}/(S_j\oplus T))
&= \{ \bar\Phi_{S_j}\circ \sigma_1\circ\varphi_j(a)\oplus
      \bar g(a) \; \vert \; a \in A_T \;\}  \cr
&=\{ \bar\Phi_{S_j}\circ \sigma_1\circ\varphi_j\circ \bar g^{-1}(a)\oplus
   a \; \vert \; a \in A_T \;\} .
\end{aligned}
$$
Comparing these two equations, we obtain
$
\sigma_2=\bar \Phi_{S_j}\circ\sigma_1\circ(\varphi_j\circ \bar g^{-1}
\circ \varphi_j^{-1})$, which implies
$[\sigma_1]=[\sigma_2]$ in $O(S_j)\setminus O(A_{S_j})/ G$.

To see that the map $\mu_j$ is well-defined on the double coset
$O(S_j)\setminus O(A_{S_j})/ G$, let us assume
$[\sigma_1]=[\sigma_2]$ in the double coset.  Then there exist
$f \in O(S_j)$ and $g \in G$ such that
$
\sigma_2=\bar f\circ\sigma_1\circ(\varphi_j\circ\bar g\circ\varphi_j^{-1})$.
Now it is clear from the above
calculations that for the $\mathbf Q$-linear extensions $f\oplus g: S_j^*\oplus
T^* \rightarrow S_j^*\oplus T^*$ we have
$
(f\oplus g)(L^{\sigma_1\circ\varphi_j})=L^{\sigma_2\circ\varphi_j}$.
Then, by Remark after Definition A.6, $L^{\sigma_1\circ\varphi_j}$ and
$L^{\sigma_2\circ\varphi_j}$ are $G$-equivalent.

These two observations show that the map $\mu_j$ descends to a well-defined
injective map $\bar\mu_j$. Since $\mu_j$ is surjective, so is
$\bar \mu_j$.
\end{proof}

\vskip0.5cm

Let $[\iota] \in \GPEj$. Then, by definition,
$\iota(T)^\perp \cong S_j$. We choose a representative
$\iota: T \hookrightarrow \Lambda$ of this $G$-equivalence class.
We also fix an isometry $\sigma: S_j \simrightarrow \iota(T)^\perp$.
Now consider the isomorphism $\sigma \oplus \iota : S_j\oplus T \simrightarrow
\iota(T)^\perp \oplus \iota(T)$ and its extension to the dual lattices.
Then one has the following commutative diagram:
$$
\begin{matrix}
\sigma \oplus \iota: & S_j^* \oplus T^* & \simrightarrow &
(\iota(T)^\perp)^* \oplus \iota(T)^* \cr
& \cup  & & \cup \cr
&L(\sigma, \iota) & \simrightarrow & \Lambda \cr
&\cup &  & \cup \cr
\sigma \oplus \iota : & S_j \oplus T & \simrightarrow &
\iota(T)^\perp \oplus \iota(T) \cr
\end{matrix}
$$
where we define the over-lattice $L(\sigma, \iota)$ of $S_j\oplus T$ by
$
L(\sigma, \iota):=(\sigma\oplus\iota)^{-1} (\Lambda)$.
Since $\Lambda$ is even integral unimodular, so is $L(\sigma,\iota)$,
i.e. $L(\sigma, \iota) \in \mathcal E_j$.

\begin{lemma} \label{lemma:lemmaLsi}
The over-lattices $L(\sigma, \iota)$ above
satisfy the following:
\begin{list}{}{
\setlength{\leftmargin}{10pt}
\setlength{\labelwidth}{6pt}
}
\item[1)] If
$\sigma: S_j \simrightarrow \iota(T)^\perp$ and
$\sigma': S_j \simrightarrow \iota(T)^\perp$, then
$
[L(\sigma,\iota)]=[L(\sigma',\iota)]
$
as an element of $\GEj$.
\item[2)] If
$\iota: T\hookrightarrow \Lambda$ and
$\iota': T\hookrightarrow \Lambda$ are $G$-equivalent primitive
embeddings $\PEj$, then
$[\iota] = [\iota']$ in $\GPEj$ and
$
[L(\sigma,\iota)]=[L(\sigma,\iota')]
$
in $\GEj$.
\end{list}
\end{lemma}

\begin{proof}
Both claims are immediate from a diagram similar to the above and
the definition
of $G$-equivalence. Details are left for readers.
\end{proof}
\begin{proposition} \label{prop:GPEbijectGEj}
\begin{list}{}{
\setlength{\leftmargin}{10pt}
\setlength{\labelwidth}{6pt}
}
\item[1)] The following map is well-defined:
$$
\xi: \GPEj   \rightarrow  \GEj , \quad
[\iota]  \mapsto  [L(\sigma,\iota)]
\label{eqn:mapxi}
$$
\item[2)] The map $\xi:  \GPEj   \rightarrow  \GEj$ is bijective.
\end{list}
\end{proposition}

\begin{proof}
The well-definedness of the map $\xi$ is immediate by Lemma A.8.
To show surjectivety, let $[K]$ be in $\GEj$. Then we have
$S_j\oplus T \subset K \subset S_j^*\oplus T^*$, and $T\hookrightarrow K$
is primitive.  Let
$\iota_K: T \hookrightarrow K$ be an inclusion.
Since $K$ is an even
unimodular indefinite lattice, there exists an isometry
$f: K \simrightarrow \Lambda$ by the theorem of Milnor (see [Se]).
Then for the embedding $\iota:=f\circ\iota_K: T\hookrightarrow
\Lambda$ and $\sigma:=f|_{S_j}: S_j \simrightarrow \iota(T)^\perp$ we have
$K=L(\sigma,\iota)$, that is $[K]=\xi([\iota])$. To show the injectivity
of $\xi$, assume $[L(\sigma,\iota)]=
[L(\sigma',\iota')]$ in $\GEj$. Then we have:
$$
\begin{matrix}
(\iota(T)^\perp)^*\oplus\iota(T)^* & \simmapleftarrow{} & S_j^*\oplus T^*
& \simmaprightarrow{} & S_j^*\oplus T^* & \simmaprightarrow{} &
(\iota'(T)^\perp)^*\oplus\iota'(T)^* \cr
\cup & & \cup & & \cup & & \cup \cr
\Lambda & \simmapleftarrow{} & L(\sigma,\iota)
& \simmaprightarrow{\exists \Phi} &
L(\sigma', \iota') & \simmaprightarrow{} & \Lambda \cr
\cup & & \cup &  & \cup & & \cup \cr
\iota(T)^\perp\oplus\iota(T)
& \simmapleftarrow{\sigma\oplus\iota} & S_j\oplus T &
\simmaprightarrow{\Phi_{S_j}\oplus\Phi_T} & S_j\oplus T
& \simmaprightarrow{\sigma'\oplus\iota'} &
\iota'(T)^\perp\oplus\iota'(T) \cr
&  & &  & & & \cr
\end{matrix}
$$
where $\Phi_{S_j}:=\Phi|_{S_j} \in O(S_j)$ and $\Phi_T:=\Phi|_T
\in G$. If we define $\Phi'=(\sigma'\oplus\iota')\circ\Phi\circ
(\sigma\oplus\iota)^{-1}: \Lambda \simrightarrow \Lambda$, we see that
$\Phi'\circ \iota=\iota'\circ\Phi_T$, which means that
$[\iota]=[\iota']$
in $\GPEj$. This completes the proof.
\end{proof}

\section{The group of Hodge isometries}

In this appendix, we shall give a proof of
the following Proposition applied in previous sections. This is probably
well-known to experts. Our proof here is similar to an argument
in [Ni2][St].

\begin{proposition}
Let $T$ be a lattice of signature
$(2, \text{rank}\, T -2)$ and $(T, \mathbf C \omega)$
be a weight two Hodge structure on $T$, i.e. $\omega$ is an element of
$T \otimes \mathbf C$ such that
$(\omega, \omega) = 0$ and $(\omega, \overline{\omega}) > 0$.
Assume the following minimality condition on $T$: if $T' \subset T$
is primitive and $\omega \in T'
\otimes \mathbf C$, then $T' = T$. Then $O_{Hodge}(T, \mathbf C \omega)$
is a cyclic group of even order, say, $2I$
such that $\varphi(2I) \vert \text{rank}\, T$. Here $\varphi(J) =
\vert (\mathbf Z/J)^{\times} \vert$ is the Euler function.
Moreover, if $g$ is a generator of $O_{Hodge}(T, \mathbf C \omega)$,
then $g(\omega) = \zeta_{2I}\omega$,
where $\zeta_{2I}$ is a primitive $2I$-th root of unity.
\end{proposition}

We proceed in five steps.

\vskip0.3cm
\noindent
{\it Step 1.} $\vert O_{Hodge}(T, \mathbf C \omega) \vert < \infty$.

\begin{proof} Put $t := \text{rank}\, T$. Let $P := \langle
\text{Re}\,\omega, \text{Im}\,\omega \rangle$. Then $P$ is a
positive definite $2$-dimensional plane by $(\omega, \omega) = 0$,
$(\omega, \overline{\omega}) > 0$.  Set $N :=
P^{\perp}$ in $T \otimes \mathbf R$. Since $\text{sgn}\,T = (2, t-2)$,
$N$ is a negative definite. Let $g \in
O_{Hodge}(T, \mathbf C \omega)$. By definition of
$O_{Hodge}(T, \mathbf C \omega)$, we have $g(P) = P$ and $g(N) = N$.
Therefore, via the assignment $g \mapsto (g \vert P, g \vert N)$, we
have an embedding $O_{Hodge}(T, \mathbf C \omega)
\subset O(P) \times O(N)$. This should be continuous. Since $P$ and $N$
are both definite, $O(P) \times O(N)$ is compact. Since
$O_{Hodge}(T, \mathbf C \omega) \subset O(T)$ is discrete, this implies
the result. \end{proof}

Let $g \in O_{Hodge}(T, \mathbf C \omega)$. Then $g(\omega) =
\alpha(g)\omega$ for some $\alpha(g) \in \mathbf C^{\times}$.
The assignment $g \mapsto \alpha(g)$ defines a group homomorphism
$\alpha : O_{Hodge}(T, \mathbf C \omega)
\rightarrow \mathbf C^{\times}$.

\vskip0.3cm
\noindent
{\it Step 2.} $\alpha$ is injective.

\begin{proof}  Let $g \in \text{Ker}\, \alpha$. Since $g$ is defined over
$\mathbf Z$ (i.e. $g$ is represented by a
matrix of integral entries with respect to integral basis of $T$),
the invarinat set $T^{g} :=
\{x \in T \vert g(x) = x\}$
is a primitive sublattice of $T$ such that $T^{g} \otimes \mathbf C
= (T \otimes \mathbf C)^{g} \ni \omega$. Thus, by
the minimality condition, we have $T^{g} = T$.
This means $g = \text{id}_T$.  \end{proof}

\vskip0.3cm
\noindent
{\it Step 3.} $O_{Hodge}(T, \mathbf C \omega) \simeq \mathbf Z/2I$ for
some $I \in \mathbf N$.

\begin{proof} By Steps 1 and 2, $O_{Hodge}(T, \mathbf C \omega)$ is
isomorphic to a finite subgroup of
the multiplicative group $\mathbf C^{\times}$ of the field $\mathbf C$.
Such groups are always cyclic. Therefore
$O_{Hodge}(T, \mathbf C \omega)$ is
a finite cyclic
group. Since the involution $-id_{T}$ is an element of $O_{Hodge}(T,
\mathbf C \omega)$, the order of
$O_{Hodge}(T, \mathbf C \omega)$ is also even.
\end{proof}

Set $O_{Hodge}(T, \mathbf C \omega) = \langle g \rangle$. Then
$\text{ord}\, (g) = 2I$.

\vskip0.3cm
\noindent
{\it Step 4.} All the eigenvalues of $g$ (on $T \otimes \mathbf C$) are
primitive $2I$-th roots of unity.

\begin{proof} Since $\text{ord}\, (g) = 2I$, all the eigenvalues of
$g$ are $2I$-th roots of unity. Since
$\alpha$ is injective, $g(\omega) = \zeta_{2I}\omega$, where $\zeta_{2I}$
is a primtive $2I$-th root of unity.
Let $d$ be a positive integer such that $d \vert 2I$ and $d \not= 2I$.
Suppose that $g$ admits a primitive $d$-th root of unity $\zeta_{d}$ as its
eigenvalue.  Then $g^{d}$ has an eigenvalue $1$. Since $g$ is defined
over $\mathbf Z$, there is $x \in T - \{0\}$ such
that $g^{d}(x) = x$. Since $g^{d}(\omega) = \zeta_{e}\omega$, where
$e = 2I/d > 1$ and $\zeta_{e}$ is a primitive $e$-th root of unity, we have
$$(x, \omega) = (g^{d}(x), g^{d}(\omega)) = (x, \zeta_{e}\omega)\, .$$
This implies $(x, \omega) = 0$, because $\zeta_{e} \not= 1$ by $e > 1$.
However, we have then $\omega \in (x)^{\perp}
\otimes \mathbf C$, a contradiction to the minimality condition of $T$.
Therefore all the eigenvalues of $g$ are primitive
$2I$-th roots of unity.
\end{proof}

\vskip0.3cm
\noindent
{\it Step 5.} $\varphi(2I) \vert t = \text{rank}\, T$.

\begin{proof}
Since $g$ is defined over $\mathbf Q$ and the
primitive $2I$-th roots of unity are mutually conjugate
over $\mathbf Q$, their multiplicities as eigenvalue of $g$ are the same,
say $m$.  Combining this with the fact
that the number of primitive $2I$-th roots of unity is
$\varphi(2I)$,  we have $m \cdot \varphi(2I) = t$.
This implies $\varphi(2I) \vert t$. \end{proof}



\begin{thebibliography}{[BKPS]}



\bibitem[BPV]{BPV} W. Barth, C. Peters, A. Van de Ven,
\textit{Compact complex surfaces }, Springer-Verlag (1984).

\bibitem[BM]{BM} T. Bridgeland, A. Maciocia,
\textit{ Complex surfaces with equivalent derived categories}, Math. Zeit.
{\bf 236} (2001) 677--697.


\bibitem[Cs]{Cs} J. W. S. Cassels,
\textit{ Rational quadratic forms},  Academic Press (1978).

\bibitem[CoS]{CoS} J. H. Conway, N. J. A. Sloane, \textit{Sphere Packings,
Lattices and Groups}, Springer-Verlag (1988).


\bibitem[Fl]{Fl} D. Flath,
\textit{Introduction to Number Theory}, Wiley-Interscience (1989).


\bibitem[HLOY1]{HLOY1}  S. Hosono, B. Lian, K. Oguiso, S.T. Yau,
\textit{Autoequivalences of derived category of a K3 surface
        and monodromy transformations }, math.AG/0201047, accepted
by J. Alg. Goem.

\bibitem[HLOY2]{HLOY2}  S. Hosono, B. Lian, K. Oguiso, S.T. Yau,
\textit{Kummer structures on a K3 surface- an old question of T. Shioda},
math.AG/0202082, accepted by Duke Math. J.


\bibitem[Ma1]{Ma1} M. Maruyama
\textit{Stable vector bundles on an algebraic surface},
Nagoya Math. J. {\bf 58} (1975) 25--68.

\bibitem[Ma2]{Ma2} M. Maruyama
\textit{Moduli of stable sheaves, II},
J. Math. Kyoto Univ. {\bf 18} (1978) 557--614.

\bibitem[MM]{MM}  R. Miranda, D. Morrison,
\textit{ The number of Embeddings of integral quadratic forms, I },
Proc. Japan Acad. Ser. A {\bf 61} (1985) 317--320.

\bibitem[Mo]{Mo}  D. Morrison,
\textit{On K3 surfaces with large Picard number },
Invent. Math. {\bf 75} (1984) 105--121.


\bibitem[Mu1]{Mu1}  S. Mukai,
\textit{Duality between $\mathbf D(X)$ and $\mathbf D(\hat X)$
        with its application to Picard sheaves },
Nagoya Math. J. {\bf 81} (1981) 101--116.


\bibitem[Mu2]{Mu2}  S. Mukai,
\textit{ Symplectic structure of the moduli space of sheaves
         on an abelian or K3 surface},
Invent. Math. {\bf 77} (1984) 101--116.


\bibitem[Mu3]{Mu3} S. Mukai,
\textit{ On the moduli space of bundles on K3 surfaces I},
in Vector bundles on algebraic varieties, Oxford Univ. Press (1987) 341--413.


\bibitem[Ni1]{Ni1}  V. V. Nikulin,
\textit{Integral symmetric bilinear forms and some of their geometric
applications }, Math. USSR Izv. {\bf 14} (1980) 103--167.

\bibitem[Ni2]{Ni2} V. V. Nikulin,
\textit{Finite groups of automorphisms of K\"ahlerian surfaces of type K3 },
Moscow Math. Sod. {\bf 38} (1980) 71--137.


\bibitem[Og1]{Og1}  K. Oguiso,
\textit{K3 surfaces via almost-primes }, math.AG/0110282,
Math. Res. Lett. {\bf 9} (2002) 47--63.


\bibitem[Og2]{Og2} K. Oguiso,
\textit{Local families of K3 surfaces and applications}, math.AG/0011258
and math.AG/0104049, accepted by J. Alg. Geom.


\bibitem[Or]{Or} D. Orlov,
\textit{Equivalences of derived categories and K3 surfaces },
math.AG/9606006,
J. Math.Sci. {\bf 84} (1997) 1361--1381.


\bibitem[Se]{Se}  J. P. Serre,
\textit{A course in arithmetic }, Springer-Verlag (1973).


\bibitem[St]{St}  H. Sterk,
\textit{Finiteness results for algebraic K3 surfaces}, Math. Zeit.
{\bf 189} (1985) 507--513.


\bibitem[Za]{Za}  D. Zagier,
\textit{Zetafunktionen und quadratische Korper : eine Einfuhrung in die
hohere Zahlentheorie },  Springer-Verlag (1981).

\end{thebibliography}
\end{document}